# TWO SUFFICIENT CONDITIONS FOR POISSON APPROXIMATIONS IN THE FERROMAGNETIC ISING MODEL

By David Coupier

Université Lille 1

A $d$-dimensional ferromagnetic Ising model on a lattice torus is considered. As the size of the lattice tends to infinity, two conditions ensuring a Poisson approximation for the distribution of the number of occurrences in the lattice of any given local configuration are suggested. The proof builds on the Stein–Chen method. The rate of the Poisson approximation and the speed of convergence to it are defined and make sense for the model. Thus, the two sufficient conditions are traduced in terms of the magnetic field and the pair potential. In particular, the Poisson approximation holds even if both potentials diverge.

**1. Introduction.** Suppose $\{I_\alpha\}_{\alpha \in \Gamma}$ is a finite family of indicator random variables, with the properties that the probabilities $\mathbb{P}(I_\alpha = 1)$ are small and that there is not too much dependence between the $I_\alpha$'s. Then, "the law of small numbers" says the sum $\sum_{\alpha \in \Gamma} I_\alpha$ should have approximately a Poisson distribution. The "birthday problem" and its variants (see Chen [5]), the theory of random graphs (see Bollobás [3] for a general reference or the famous paper of Erdös and Rényi [11]) and the study of words in long DNA sequences (see, e.g., Schbath [18]) are examples in which a law of small numbers takes place. As the situation studied in this paper, these examples can be viewed as problems of increasing size (i.e., the cardinality of $\Gamma$ tends to infinity) in which the sum $\sum_{\alpha \in \Gamma} I_\alpha$ has a Poisson limit. Two methods are often used for proving Poisson approximations; the moment method (see [3], page 25) and the Stein–Chen method (see Arratia, Goldstein and Gordon [1], Barbour, Holst and Janson [2] for a very complete reference, or [5] for the original paper of Chen). The second one offers two main advantages. Only the first two moments need to be computed and a bound of the rate of convergence is obtained. However, the Stein–Chen method requires to restrict our attention to variables which satisfy the FKG inequality [14]. This is the case of spins of a ferromagnetic Ising model.

Let us consider a lattice graph in dimension $d \geq 1$, with periodic boundary conditions (lattice torus). The vertex set is $V_n = \{0, \ldots, n-1\}^d$. The integer $n$ will be called the *size* of the lattice. The edge set, denoted by $E_n$, will be specified by defining the set of neighbors $\mathcal{V}(x)$ of a given vertex $x$:

$$\mathcal{V}(x) = \{y \neq x \in V_n, \|y - x\|_p \leq \rho\}, \tag{1}$$

where the substraction is taken componentwise modulo $n$, $\|\cdot\|_p$ stands for the $L_p$ norm in $\mathbb{R}^d$ ($1 \leq p \leq \infty$), and $\rho$ is a fixed parameter. For instance, the square lattice is obtained for $p = \rho = 1$. Replacing the $L_1$ norm by the $L_\infty$ norm adds the diagonals. From now on, all operations on vertices will be understood modulo $n$. In particular, each vertex of the lattice has the same number of neighbors; we denote by $\mathcal{V}$ this number.

A *configuration* is a mapping from the vertex set $V_n$ to the state space $\{-1, +1\}$. Their set is denoted by $\mathcal{X}_n = \{-1, +1\}^{V_n}$ and called the *configuration set*. The Ising model is classically defined as follows (see, e.g., Georgii [16] and Malyshev and Minlos [17]).

DEFINITION 1.1. Let $G_n = (V_n, E_n)$ be the undirected graph structure with finite vertex set $V_n$ and edge set $E_n$. Let $a$ and $b$ be two reals. The *Ising model* with parameters $a$ and $b$ is the probability measure $\mu_{a,b}$ on $\mathcal{X}_n = \{-1, +1\}^{V_n}$ defined by: $\forall \sigma \in \mathcal{X}_n$,

$$\mu_{a,b}(\sigma) = \frac{1}{Z_{a,b}} \exp\left( a \sum_{x \in V_n} \sigma(x) + b \sum_{\{x,y\} \in E_n} \sigma(x)\sigma(y) \right), \tag{2}$$

where the normalizing constant $Z_{a,b}$ is such that $\sum_{\sigma \in \mathcal{X}_n} \mu_{a,b}(\sigma) = 1$.

Following the definition of Malyshev and Minlos ([17], page 2), the measure $\mu_{a,b}$ defined above is a Gibbs measure associated to potentials $a$ and $b$. Expectations relative to $\mu_{a,b}$ will be denoted by $\mathbb{E}_{a,b}$.

In the classical presentation of statistical physics, the elements of $\mathcal{X}_n$ are spin configurations; each vertex of $V_n$ is an atom whose spin is either positive or negative. Here, we shall simply talk about positive or negative vertices instead of positive or negative spins and we shall merely denote by $+$ and $-$ the states $+1$ and $-1$. The parameters $a$ and $b$ are, respectively, the *magnetic field* and the *pair potential*. The model remaining unchanged by swapping positive and negative vertices and replacing $a$ by $-a$, we chose to study only negative values of the magnetic field $a$. Moreover, throughout the paper, the pair potential $b$ will be supposed nonnegative.

Many laws of small numbers have been already proved for the Ising model. They are obtained in two different contexts; the low temperature case (i.e., $b$ large enough) and the large magnetic field case. Chazottes and Redig [4] have obtained a Poisson law for the number of occurrences of large cylindrical events. Their result concerns the low temperature case and is based on an argument of disagreement percolation. Fernández, Ferrari and Garcia [12] have established the asymptotic Poisson distribution of contours in the nearest-neighbor Ising model at low temperature and zero magnetic field. In the same context, Ferrari and Picco [13] have applied the Stein–Chen method to obtain bounds on the total variation distance between the law of large contours and a Poisson process. Ganesh et al. [15] have studied the Ising model for positive and fixed values of $b$. Provided the magnetic field $a$ tends to infinity, they proved that the distribution of the number of negative vertices is approximately Poisson. In Coupier [6], we have generalized this result to (deterministic) objects more elaborated than a single vertex, called the *local configurations*; but the parameter $b$ was still fixed.

Our goal is to extend these Poisson approximations to values of the pair potential $b$ not necessarily bounded. In particular, we establish a law of small numbers in the case where both magnetic field $a$ and pair potential $b$ diverge. Moreover, compared with the previous articles, the results of this paper offer two advantages. First, the roles of the rate of the Poisson approximation and the speed of convergence to it are conceptually clearer. Second, the method of the proof is more elementary and simple than those of the previous references; it does not require the need of auxiliary birth-and-death processes as in [12] and [13], or of activity expansions as in [15].

Pattern recognition and image denoising problems motivate the study realized in this paper. Indeed, Desolneux, Moisan and Morel (see [9] and [10]) describe and detect the geometric structures of an image through the notion of "meaningful events." Their method is based on the link between the perception threshold of a given visual structure and its probability of appearance in a random image. In Coupier, Desolneux and Ycart [7], a Poisson approximation for the probability of appearance of any local pattern is computed in the case of independent pixels. This result leads to a denoising algorithm for gray-level images. However, the structures created by the noise are generally more irregular (with a large perimeter) than the natural ones of an image. This remark justifies our interest for Poisson approximations in probabilistic models with dependent pixels, as the Ising model.

We are interested in the occurrences in the graph $G_n$ of a fixed local configuration $\eta$ (see Section 2 for a precise definition and Figure 1 for an example). Such a configuration is called "local" in the sense that the vertex set on which it is defined is fixed and does not depend on $n$. Its number of occurrences in $G_n$ is denoted by $X_n(\eta)$. A local configuration $\eta$ is determined by its set of positive vertices $V_+(\eta)$ whose cardinality and perimeter are respectively denoted by $k(\eta)$ and $\gamma(\eta)$. As the size $n$ of the lattice tends to infinity, the potentials $a = a(n)$ and $b = b(n)$ depend on $n$. The case where $a(n)$ tends to $-\infty$ corresponds to rare positive vertices among a majority of negative ones. In order to simplify formulas, the Gibbs measure $\mu_{a(n),b(n)}$ will be merely denoted by $\mu_{a,b}$.

A natural idea consists in regarding both potentials $a(n)$ and $b(n)$ through the same quantity: the *weight* of the local configuration $\eta$

$$W_n(\eta) = \exp(2a(n)k(\eta) - 2b(n)\gamma(\eta)).$$

Actually, the weight $W_n(\eta)$ represents the probabilistic cost associated to a given occurrence of $\eta$ and, consequently, the product $n^d W_n(\eta)$ (where $n^d$ is the cardinality of $V_n$) corresponds to the expected number of occurrences of $\eta$ in the whole graph $G_n$. Hence, in order to obtain a Poisson approximation for the random variable $X_n(\eta)$ it is needed that the weight $W_n(\eta)$ tends to 0 at the rate $n^{-d}$. Therefore, throughout this paper, the potentials $a(n)$ and $b(n)$ will satisfy the homogeneity relation

$$(3) \qquad n^d W_n(\eta) = \lambda,$$

where $\lambda$ is a positive constant. Two other parameters naturally take place in our study: the probability gap $\Delta_n(\eta)$ and the maximality probability $\Theta_n(\eta)$. When the first one tends to 0, the vertices surrounding a given occurrence of $\eta$ are all negative. A null limit for the second one means if somewhere in the graph the vertices corresponding to the set $V_+(\eta)$ are positive, then locally they are the only ones. The fact that the probability gap $\Delta_n(\eta)$ and the maximality probability $\Theta_n(\eta)$ tend to 0 is denoted by (H1). We denote by (H2) the following (interpreted) hypothesis: "the probability for a given configuration $\zeta \in \{-,+\}^V$, $V \subset V_n$, to occur is bounded by its weight $W_n(\zeta)$." Our main result (Theorem 3.1) describes the asymptotic behavior of the number of occurrences of $\eta$ in the lattice. Added to the necessary condition (3), the two hypotheses (H1) and (H2) imply the total variation distance between the distribution of the random variable $X_n(\eta)$ and the Poisson distribution with parameter $\lambda$ satisfies

$$d_{\mathrm{TV}}(\mathcal{L}(X_n(\eta)), \mathcal{P}(\lambda)) = \mathcal{O}(\max\{\Delta_n(\eta), \Theta_n(\eta)\}).$$

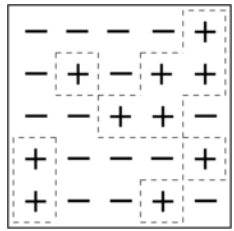

FIG. 1. *A local configuration $\eta$ with $k(\eta) = |V_+(\eta)| = 10$ positive vertices and a perimeter $\gamma(\eta)$ equal to 58, in dimension $d = 2$ and on a ball of radius $r = 2$ (with $\rho = 1$ and relative to the $L_\infty$ norm).*

The proof is based on the Stein–Chen method. Lemma 4.1 reduces the proof to a sum of *increasing* random indicators $\overline{X}_n(\eta)$. Since the Gibbs measure $\mu_{a,b}$ defined in (2) satisfies the FKG inequality (because the pair potential $b$ is positive), this method is applied to the random variable $\overline{X}_n(\eta)$ and produces Lemma 4.5. Then, bounds on the first two moments of $\overline{X}_n(\eta)$ (Lemmas 4.2 and 4.6) allow to conclude.

Finally, Theorem 3.1 is completed by a second result (Proposition 3.2) which traduces the two hypotheses (H1) and (H2) in terms of magnetic field $a(n)$ and pair potential $b(n)$.

The paper is organized as follows. The notion of local configuration $\eta$ is defined in Section 2. Its number of positive vertices $k(\eta)$, its perimeter $\gamma(\eta)$ and its weight $W_n(\eta)$ are also introduced. The use of the *local energy* of $\eta$ (Definition 2.1) allows to write the conditional probability for $\eta$ to occur in the graph in terms of weights (Lemma 2.2). Lemma 2.4 is the main result of Section 2. It provides lower and upper bounds for this conditional probability on which the estimates of the first two moments of $\overline{X}_n(\eta)$ are based. In Section 3, Theorem 3.1 and hypotheses (H1) and (H2) are introduced and discussed. Thus, Proposition 3.2 is illustrated by two examples. Finally, Section 4, in which the sum $\overline{X}_n(\eta)$ is defined, is devoted to the proof of Theorem 3.1.

**2. Conditional probability of a local configuration.** Let us start with some notation and definitions. Given $\zeta \in \mathcal{X}_n = \{-,+\}^{V_n}$ and $V \subset V_n$, we denote by $\zeta_V$ the natural projection of $\zeta$ over $\{-,+\}^V$. If $U$ and $V$ are two disjoint subsets of $V_n$, then $\zeta_U \zeta'_V$ is the configuration on $U \cup V$ which is equal to $\zeta$ on $U$ and $\zeta'$ on $V$. Let us denote by $\delta V$ the neighborhood of $V$ [corresponding to (1)]:

$$\delta V = \{y \in V_n \setminus V, \exists x \in V, \{x,y\} \in E_n\},$$

and by $\overline{V}$ the union of the two disjoint sets $V$ and $\delta V$. Moreover, $|V|$ denotes the cardinality of $V$ and $\mathcal{F}(V)$ the $\sigma$-algebra generated by the configurations of $\{-,+\}^V$.

As usual, the graph distance dist is defined as the minimal length of a path between two vertices. We shall denote by $B(x,r)$ the ball of center $x$ and radius $r$:

$$B(x,r) = \{y \in V_n; \mathrm{dist}(x,y) \leq r\}.$$

In the case of balls, $\overline{B(x,r)} = B(x,r+1)$. In order to avoid unpleasant situations, like self-overlapping balls, we will always assume that $n > 2\rho r$. If $n$ and $n'$ are both larger than $2\rho r$, the balls $B(x,r)$ in $G_n$ and $G_{n'}$ are isomorphic. Two properties of the balls $B(x,r)$ will be crucial in what follows. The first one is that two balls with the same radius are translates of each other:

$$B(x+y, r) = y + B(x, r).$$

The second one is that for $n > 2\rho r$, the cardinality of $B(x,r)$ depends only on $r$ and neither on $x$ nor on $n$: it will be denoted by $\beta(r)$. The same is true for the number of edges $\{y,z\} \in E_n$ with $y, z \in B(x,r)$, which will be denoted by $\alpha(r)$.

Let $r$ be a positive integer, and consider a fixed ball with radius $r$, say $B(0,r)$. We denote by $\mathcal{C}_r = \{-,+\}^{B(0,r)}$ the set of configurations on that ball. Elements of $\mathcal{C}_r$ will be called *local configurations of radius $r$*, or merely *local configurations* whenever the radius $r$ will be fixed. Of course, there exists only a finite number of such configurations (precisely $2^{\beta(r)}$). See Figure 1 for an example. Throughout this paper, the radius $r$ will be constant, that is, it will not depend on the size $n$. Hence, defining local configurations on balls of radius $r$ serves only to ensure that studied objects are "local." In what follows, $\eta$, $\eta'$ will denote local configurations of radius $r$.

A local configuration $\eta \in \mathcal{C}_r$ is determined by its subset $V_+(\eta) \subset B(0,r)$ of positive vertices:

$$V_+(\eta) = \{x \in B(0,r), \eta(x) = +\}.$$

The cardinality of this set will be denoted by $k(\eta)$ and its complementary set in $B(0,r)$, that is, the set of negative vertices of $\eta$, by $V_-(\eta)$. Moreover, the geometry (in the sense of the graph structure) of the set $V_+(\eta)$ needs to be described. Let us define the *perimeter* $\gamma(\eta)$ of the local configuration $\eta$ by the formula

$$\gamma(\eta) = \mathcal{V}|V_+(\eta)| - 2|\{\{x,y\} \in V_+(\eta) \times V_+(\eta), \{x,y\} \in E_n\}|,$$

where $\mathcal{V}$ is the number of neighbors of a vertex. In other words, $\gamma(\eta)$ counts the pairs of neighboring vertices $x$ and $y$ of $B(0,r)$ having opposite spins (under $\eta$) and those such that $x \in B(0,r)$, $y \in \delta B(0,r)$ and $\eta(x) = +$. Finally, we denote by $W_n(\eta)$ and call the *weight* of the local configuration $\eta$ the following quantity:

$$W_n(\eta) = \exp(2a(n)k(\eta) - 2b(n)\gamma(\eta)).$$

Since $a(n) < 0$ and $b(n) \geq 0$, the weight $W_n(\eta)$ satisfies $0 < W_n(\eta) \leq 1$. That of the local configuration having only negative vertices, denoted by $\eta^0$ and called the *null local configuration*, is equal to 1. If $\eta \neq \eta^0$, then $k(\eta) \geq 1$ and $\gamma(\eta) \geq \mathcal{V}$. It follows that

$$W_n(\eta) \leq \exp(2a(n) - 2b(n)\mathcal{V}).$$

Actually, the weight $W_n(\eta)$ represents the probabilistic cost associated to the presence of $\eta$ on a given ball. Lemma 2.4 will give a rigorous sense to this idea.

Remark that the notation $k(\cdot)$, $\gamma(\cdot)$ and $W_n(\cdot)$ can be naturally extended to any configuration $\zeta \in \{-,+\}^V$, $V \subset V_n$. Furthermore, the configuration of $\{-,+\}^V$ having no positive vertex will be said to be *null* and denoted by $\zeta^0$.

Let $\eta \in \mathcal{C}_r$. For each vertex $x \in V_n$, denote by $\eta_x$ the translation of $\eta$ onto the ball $B(x,r)$ (up to periodic boundary conditions):

$$\forall y \in V_n \qquad \text{dist}(0,y) \leq r \implies \eta_x(x+y) = \eta(y).$$

In particular, $V_+(\eta_x) = x + V_+(\eta)$. So, $\eta$ and $\eta_x$ have the same number of positive vertices and the same perimeter. So do their weights. Let us denote by $I_x^\eta$ the indicator function defined on $\mathcal{X}_n$ as follows: $I_x^\eta(\sigma)$ is 1 if the restriction of the configuration $\sigma \in \mathcal{X}_n$ to the ball $B(x,r)$ is $\eta_x$ and 0 otherwise. Finally, let us define the random variable $X_n(\eta)$ which counts the number of copies of the local configuration $\eta$ in $G_n$:

$$X_n(\eta) = \sum_{x \in V_n} I_x^\eta.$$

Due to periodicity, this sum consists of $n^d$ indicator functions $I_x^\eta$, which have the same distribution. Hence, in order to understand the behavior of the random variable $X_n(\eta)$ it suffices to study that of one of the random indicators $I_x^\eta$. This is the reason why the rest of this section is devoted to the study of the probability for $\eta$ (or $\eta_x$) to occur on a ball $B(x,r)$ knowing what happens on its neighborhood:

(4) $$\mu_{a,b}(I_x^\eta = 1 | \sigma_{\delta B(x,r)}).$$

Besides, let us remark that the expected number of occurrences of the local configuration $\eta$ in the graph can be expressed as the following expectation:

(5) $$\begin{aligned}\mathbb{E}_{a,b}[X_n(\eta)] &= \mathbb{E}_{a,b}[n^d I_x^\eta] \\ &= \mathbb{E}_{a,b}[n^d \mu_{a,b}(I_x^\eta = 1 | \mathcal{F}(\delta B))],\end{aligned}$$

where $x$ is a given vertex and $B$ denotes merely the ball $B(x,r)$. Here, $\mu_{a,b}(I_x^\eta = 1 | \mathcal{F}(\delta B))$ represents a $\mathcal{F}(\delta B)$-measurable random variable and, for any configuration $\sigma$ of $\{-,+\}^{\delta B}$, $\mu_{a,b}(I_x^\eta = 1 | \mathcal{F}(\delta B))(\sigma)$ is equal to the conditional probability $\mu_{a,b}(I_x^\eta = 1 | \sigma)$.

Lemma 2.4 provides inequalities between the conditional probability (4) and the weight $W_n(\eta)$. A way to link these two quantities consists in using the *local energy* of $\eta$. Let us start with the following definition.

DEFINITION 2.1. Let $V \subset V_n$ and $\zeta \in \{-,+\}^{\overline{V}}$. The *local energy* $H^V(\zeta)$ of the configuration $\zeta$ on the set $V$ is defined by

$$H^V(\zeta) = a(n) \sum_{y \in V} \zeta(y) + b(n) \sum_{\substack{\{y,z\} \in E_n \\ (y \in V) \vee (z \in V)}} \zeta(y)\zeta(z),$$

where $(y \in V) \vee (z \in V)$ means at least one of the two vertices $y$ and $z$ belongs to $V$ (the other might belong to its neighborhood $\delta V$).

Let us fix a vertex $x$ and denote by $B$ the ball $B(x,r)$. On the one hand, remark that for any $\sigma \in \{-,+\}^{\delta B}$ the local energy $H^B(\eta_x \sigma)$ on $B$ of the configuration which is equal to $\eta_x$ on $B$ and $\sigma$ on $\delta B$ gives an explicit formula of the conditional probability $\mu_{a,b}(I_x^\eta = 1|\sigma)$:

(6) $$\mu_{a,b}(I_x^\eta = 1|\sigma) = \frac{\exp(H^B(\eta_x \sigma))}{\sum_{\eta' \in \mathcal{C}_r} \exp(H^B(\eta'_x \sigma))}.$$

On the other hand, the exponential of the local energy of $\eta$ can be expressed as a function of the weight $W_n(\eta)$. Indeed, $H^B(\eta_x \sigma)$ is equal to

(7) $$a(n)(2k(\eta) - \beta(r)) + b(n)\left( \sum_{\substack{\{y,z\} \in E_n \\ y,z \in B}} \eta_x(y)\eta_x(z) + \sum_{\substack{\{y,z\} \in E_n \\ y \in B, z \in \delta B}} \eta_x(y)\sigma(z) \right).$$

Let $\{y,z\}$ be an edge such that $y \in B$ and $\eta_x(y) = +$. Assume $z \in B$. If $\eta_x(z) = -$, the edge $\{y,z\}$ is counted by the perimeter $\gamma(\eta)$, whereas if $\eta_x(z) = +$, it is not. In the case where $z \notin B$ (i.e., $z \in \delta B$), the edge $\{y,z\}$ is consistently counted by $\gamma(\eta)$. In other words, the following equality holds:

(8) $$2\gamma(\eta) = \alpha(r) - \sum_{\substack{\{y,z\} \in E_n \\ y,z \in B}} \eta_x(y)\eta_x(z) + \sum_{\substack{\{y,z\} \in E_n \\ y \in B, z \in \delta B}} (1 + \eta_x(y)).$$

So, the perimeter $\gamma(\eta)$ can be inserted in the expression of $H^B(\eta_x \sigma)$ given by (7) and we get

(9) $$\begin{aligned}\exp(H^B(\eta_x \sigma)) \\ = W_n(\eta) \exp\Bigg( -a(n)\beta(r) \\ + b(n)\bigg( \alpha(r) + \sum_{\substack{\{y,z\} \in E_n \\ y \in B, z \in \delta B}} 1 + \eta_x(y)(1 + \sigma(z)) \bigg) \Bigg).\end{aligned}$$

Remark that the previous identity remains valid replacing the ball $B$ with any given set $V$ and the configuration $\eta_x \sigma$ with any given $\zeta \in \{-,+\}^{\overline{V}}$.

Using the previous relations, an alternate expression of the conditional probability $\mu_{a,b}(I_x^\eta = 1|\sigma)$ is obtained.

LEMMA 2.2. *Let $x \in V_n$ and $B = B(x,r)$. Then, for any $\eta \in \mathcal{C}_r$ and $\sigma \in \{-,+\}^{V_n \setminus B}$,*

$$\mu_{a,b}(I_x^\eta = 1|\sigma_{\delta B}) = \frac{W_n(\eta_x \sigma_{\delta B})}{\sum_{\eta' \in \mathcal{C}_r} W_n(\eta'_x \sigma_{\delta B})}. \tag{10}$$

PROOF. Let $\eta$ be a local configuration of radius $r$ and $\sigma \in \{-,+\}^{V_n \setminus B}$. Let us denote by $D$ the set $\delta B \cup \delta \overline{B}$: a vertex $y$ belonging to $D$ satisfies $r < \text{dist}(x,y) \leq r+2$. Adding terms which only depend on the configuration $\sigma$ to the right-hand side of (6) does not change the equality:

$$\mu_{a,b}(I_x^\eta = 1|\sigma_{\delta B}) = \frac{\exp(H^{\overline{B}}(\eta_x \sigma_D))}{\sum_{\eta' \in \mathcal{C}_r} \exp(H^{\overline{B}}(\eta'_x \sigma_D))}. \tag{11}$$

Let $\eta' \in \mathcal{C}_r$. Applying relation (9) to the configuration $\eta'_x \sigma_D$, we get

$$\exp(H^{\overline{B}}(\eta'_x \sigma_D))$$
$$= W_n(\eta'_x \sigma_{\delta B}) \exp\Bigg(-a(n)\beta(r+1)$$
$$+ b(n)\bigg(\alpha(r+1) + \sum_{\substack{\{y,z\} \in E_n \\ y \in \overline{B}, z \in \delta \overline{B}}} 1 + \sigma_D(y)(1 + \sigma_D(z))\bigg)\Bigg).$$

Finally, it suffices to remark that the latter exponential does not depend on $\eta'_x$ and to simplify the right-hand side of (11). □

Let $V$ and $V'$ be two disjoint subsets of vertices. The following relations:

$$k(\zeta \zeta') = k(\zeta) + k(\zeta') \quad \text{and} \quad \gamma(\zeta \zeta') \leq \gamma(\zeta) + \gamma(\zeta')$$

are true whatever the configurations $\zeta \in \{-,+\}^V$ and $\zeta' \in \{-,+\}^{V'}$. The *connection* between $\zeta \in \{-,+\}^V$ and $\zeta' \in \{-,+\}^{V'}$, denoted by $\text{conn}(\zeta, \zeta')$, is defined by

$$\text{conn}(\zeta, \zeta') = |\{\{y,z\} \in E_n, y \in V, z \in V' \text{ and } \zeta(y) = \zeta'(z) = +\}|.$$

This quantity allows to link the perimeters of the configurations $\zeta \zeta'$, $\zeta$ and $\zeta'$ together:

$$\gamma(\zeta \zeta') + 2\text{conn}(\zeta, \zeta') = \gamma(\zeta) + \gamma(\zeta'),$$

and as a consequence their weights:

$$W_n(\zeta \zeta') \exp(-4b(n) \text{conn}(\zeta, \zeta')) = W_n(\zeta) W_n(\zeta'). \tag{12}$$

In particular, if the connection $\mathrm{conn}(\zeta,\zeta')$ is null, then the weight $W_n(\zeta\zeta')$ is equal to the product $W_n(\zeta)W_n(\zeta')$. This is the case when $\overline{V}\cap V'=\varnothing$.

Let $B=B(x,r)$ and $\eta$ be a local configuration of radius $r$. Assume the set of positive vertices of $\eta$ satisfies $V_+(\eta)\subset B(0,r-1)$; such a local configuration is said to be *clean*. Then, the connection $\mathrm{conn}(\eta_x,\sigma)$ is null and $W_n(\eta_x\sigma)$ is equal to the product $W_n(\eta)W_n(\sigma)$. This is the case also for the null configuration $\sigma^0$; $W_n(\eta_x\sigma^0)=W_n(\eta)W_n(\sigma^0)=W_n(\eta)$.

Let us define the *probability gap* of $\eta$ and denote by $\Delta_n(\eta)$ the quantity below:

$$\Delta_n(\eta) = \max_{\substack{\sigma\in\delta B \\ \sigma\neq\sigma^0}} \frac{W_n(\eta_x\sigma)}{W_n(\eta)}.$$

Since the weights of a configuration and its translates are equal, the probability gap $\Delta_n(\eta)$ does not depend on the ball $B=B(x,r)$ (nor on $\delta B$). It only depends on the local configuration $\eta$ and on the potentials $a(n)$ and $b(n)$. Moreover, $\Delta_n(\eta)$ satisfies the following inequalities:

LEMMA 2.3. *Let $\eta\in\mathcal{C}_r$. Then, the probability gap $\Delta_n(\eta)$ satisfies*

(13) $$\exp(2a(n)-2b(n)\mathcal{V}) \leq \Delta_n(\eta) \leq \exp(2a(n)).$$

*In particular, for all $\eta'\in\mathcal{C}_r$, $\eta'\neq\eta^0$:*

$$W_n(\eta') \leq \Delta_n(\eta).$$

PROOF. Let $\eta\in\mathcal{C}_r$ and as usual denote by $B$ the ball $B(x,r)$. Since $a(n)<0$ and $b(n)\geq 0$, the weight $W_n(\eta_x\sigma)$ is a decreasing function of parameters $k(\eta_x\sigma)$ and $\gamma(\eta_x\sigma)$. If $\sigma$ and $\sigma'$ are two configurations of $\{-,+\}^{\delta B}$, different from $\sigma^0$, such that $V_+(\sigma)\subset V_+(\sigma')$ and $k(\sigma)=k(\sigma')-1$, then the perimeter $\gamma(\eta_x\sigma)$ is not necessarily smaller than $\gamma(\eta_x\sigma')$. However, the following statement based on the convexity of balls on which the local configurations are defined is true: for all configuration $\sigma\in\{-,+\}^{\delta B}$, $\sigma\neq\sigma^0$, there exists a configuration $\sigma_{\min}(\sigma)\in\{-,+\}^{\delta B}$ satisfying

$$k(\sigma_{\min}(\sigma))=1, \qquad V_+(\sigma_{\min}(\sigma))\subset V_+(\sigma) \quad \text{and} \quad \gamma(\eta_x\sigma_{\min}(\sigma))\leq \gamma(\eta_x\sigma).$$

As a consequence, the maximum of $W_n(\eta_x\sigma)$ is obtained among the configurations $\sigma\in\{-,+\}^{\delta B}$ having only one positive vertex; $k(\sigma)=1$ and $\gamma(\sigma)=\mathcal{V}$. For such a configuration $\sigma$,

$$\frac{W_n(\eta_x\sigma)}{W_n(\eta)} = \exp(2a(n)-2b(n)(\gamma(\eta_x\sigma)-\gamma(\eta)))$$
$$= \exp(2a(n)-2b(n)(\mathcal{V}-2\,\mathrm{conn}(\eta_x,\sigma)))$$
$$\geq \exp(2a(n)-2b(n)\mathcal{V}).$$

The lower bound of (13) follows. The upper bound is also a consequence of the convexity of the ball $B$: the perimeter $\gamma(\eta_x \sigma)$, $\sigma \neq \sigma^0$, is necessarily as large as $\gamma(\eta_x \sigma^0)$. In other words, the difference $\gamma(\eta_x \sigma) - \gamma(\eta)$ is positive.

Finally, we have already seen that the weight of a configuration $\eta' \in \mathcal{C}_r$, different from $\eta^0$, is bounded by $\exp(2a(n) - 2b(n)\mathcal{V})$. So, it is bounded by $\Delta_n(\eta)$. □

For example, the probability gap of a clean local configuration is equal to the lower bound given by (13), that is, $\exp(2a(n) - 2b(n)\mathcal{V})$. Indeed, if $\eta$ is clean, the connection $\mathrm{conn}(\eta_x, \cdot)$ is null and the equalities of the previous proof allow us to conclude.

Lemma 2.4 is the main result of this section. It provides lower and upper bounds for the conditional probability $\mu_{a,b}(I_x^\eta = 1|\sigma)$ depending on the weight $W_n(\eta)$ and the probability gap $\Delta_n(\eta)$:

LEMMA 2.4. *Let $\eta$ be a local configuration of radius $r$ and let $x$ be a vertex. Let us denote by $B$ the ball $B(x,r)$. Then, for all configuration $\sigma \in \{-,+\}^{\delta B}$,*

$$\mu_{a,b}(I_x^\eta = 1|\sigma) \geq W_n(\eta) \mathbb{1}_{\sigma=\sigma^0}(1 - |\mathcal{C}_r|\Delta_n(\eta)) \tag{14}$$

*and*

$$\mu_{a,b}(I_x^\eta = 1|\sigma) \leq W_n(\eta)\left(1 + \mathbb{1}_{\sigma \neq \sigma^0} \frac{\Delta_n(\eta)}{W_n(\sigma)}\right). \tag{15}$$

Under the hypotheses of Theorem 3.1, the probability gap $\Delta_n(\eta)$ tends to 0 and the probability of the indicator $\mathbb{1}_{\sigma=\sigma^0}$ tends to 1. Henceforth, relations (14) and (15) mean the weight $W_n(\eta)$ constitutes a good estimate for the probability $\mu_{a,b}(I_x^\eta = 1)$.

Before proving Lemma 2.4, let us explain why the results of this section do not hold when the pair potential is negative. If $b(n) < 0$, the condition $a(n) - b(n)\mathcal{V} \leq 0$ (which ensures that the weight of a local configuration is smaller than 1) implies the weight $W_n(\eta)$ and the probability gap $\Delta_n(\eta)$ are both bounded by $\exp(2a(n) - 2b(n)\mathcal{V})$. In particular, we lose the inequality $W_n(\eta') \leq \Delta_n(\eta)$ of Lemma 2.3 and therefore the lower bound of Lemma 2.4.

PROOF OF LEMMA 2.4. Let $\eta \in \mathcal{C}_r$ and $x \in V_n$. Let $\sigma$ be a configuration of $\{-,+\}^{\delta B}$. If $\sigma = \sigma^0$, then Lemma 2.2 implies that

$$\mu_{a,b}(I_x^\eta = 1|\sigma^0) = \frac{W_n(\eta_x \sigma^0)}{\sum_{\eta' \in \mathcal{C}_r} W_n(\eta'_x \sigma^0)}$$

$$= \frac{W_n(\eta)}{\sum_{\eta' \in \mathcal{C}_r} W_n(\eta')}$$

$$= \frac{W_n(\eta)}{1 + \sum_{\eta' \neq \eta^0} W_n(\eta')}$$

$$\leq W_n(\eta),$$

whereas, if $\sigma \neq \sigma^0$,

$$\mu_{a,b}(I_x^\eta = 1|\sigma) = W_n(\eta) \frac{W_n(\eta_x \sigma)}{W_n(\eta)} \frac{1}{\sum_{\eta' \in \mathcal{C}_r} W_n(\eta'_x \sigma)}$$

$$\leq W_n(\eta) \Delta_n(\eta) \frac{1}{W_n(\eta_x^0 \sigma)}$$

$$\leq W_n(\eta) \frac{\Delta_n(\eta)}{W_n(\sigma)}.$$

The upper bound (15) is deduced from the above inequalities. The lower bound is also based on Lemma 2.2. Indeed,

$$\mu_{a,b}(I_x^\eta = 1|\sigma) \geq \mathbb{1}_{\sigma=\sigma^0} \frac{W_n(\eta_x \sigma^0)}{\sum_{\eta' \in \mathcal{C}_r} W_n(\eta'_x \sigma^0)}$$

$$\geq \mathbb{1}_{\sigma=\sigma^0} \frac{W_n(\eta)}{1 + \sum_{\eta' \neq \eta^0} W_n(\eta')}.$$

Now, the weight of a local configuration $\eta' \neq \eta^0$ satisfies $W_n(\eta') \leq \Delta_n(\eta)$ (see Lemma 2.3). Hence, we can write

$$\mu_{a,b}(I_x^\eta = 1|\sigma) \geq \mathbb{1}_{\sigma=\sigma^0} \frac{W_n(\eta)}{1 + |\mathcal{C}_r|\Delta_n(\eta)}.$$

Finally, the lower bound (14) is obtained using the classical inequality

$$\forall u > -1 \qquad \frac{1}{1+u} \geq 1 - u. \qquad \square$$

**3. Poisson approximation.** This section is devoted to the main result of the paper, that is, Theorem 3.1. Let us introduce two hypotheses (H1) and (H2) that will be sufficient to prove Theorem 3.1.

For any given local configuration $\eta$, let us define the subset $\mathcal{C}_r(\eta)$ of $\mathcal{C}_r$ by

$$\mathcal{C}_r(\eta) = \{\eta' \in \mathcal{C}_r, V_+(\eta') \supset V_+(\eta)\}.$$

Each element of $\mathcal{C}_r^*(\eta) = \mathcal{C}_r(\eta) \setminus \{\eta\}$ has at least $k(\eta) + 1$ positive vertices. Thus, we denote by $\Theta_n(\eta)$ and we call the *maximality probability* the following quantity:

$$\Theta_n(\eta) = \max_{\eta' \in \mathcal{C}_r^*(\eta)} \frac{W_n(\eta')}{W_n(\eta)}.$$

The necessary condition $n^d W_n(\eta) = \lambda$, for some constant $\lambda > 0$, forces the weight $W_n(\eta)$ to tend to 0 and therefore the local configuration $\eta$ to be different from $\eta^0$. Thus, using the inequality $\gamma(\eta) \leq \mathcal{V}k(\eta)$, we deduce that the probabilistic cost associated to a single positive vertex, that is,

$$\exp(2a(n) - 2b(n)\mathcal{V}),$$

is bounded by $W_n(\eta)^{1/k(\eta)}$ and tends to 0 as $n$ tends to infinity. Actually, Theorem 3.1 requires stronger hypotheses. First, if the vertices corresponding to $x + V_+(\eta)$ are positive, then the other ones belonging to $B(x,r)$ and those belonging to the neighborhood $\delta B(x,r)$ must be negative. This fact results in

$$(\text{H1}): \lim_{n \to +\infty} \max\{\Delta_n(\eta), \Theta_n(\eta)\} = 0.$$

For $V \subset V_n$ and $\zeta \in \{-,+\}^V$, the indicator function $I_V^\zeta$ is defined as follows: $I_V^\zeta(\sigma)$ is 1 if the restriction of the configuration $\sigma \in \mathcal{X}_n$ to $V$ is $\zeta$ and 0 otherwise. Theorem 3.1 needs a control of the probability of the event $I_V^\zeta = 1$. Let (H2) be the following hypothesis:

$$(\text{H2}): \begin{array}{c} \forall V \subset V_n, \exists N(V) \in \mathbb{N}, \exists C = C(|V|) > 0, \forall \zeta \in \{-,+\}^V, \\ \forall n \geq N(V), \mu_{a,b}(I_V^\zeta = 1) \leq CW_n(\zeta). \end{array}$$

Here, $C = C(|V|)$ means that $C$ depends on the set $V$ only through its cardinality.

Before introducing our main result, let us recall some classical notations. If $\mu$ and $\nu$ are two probability distributions, the *total variation distance* between $\mu$ and $\nu$ is

$$d_{\text{TV}}(\mu, \nu) = \sup_A |\mu(A) - \nu(A)|,$$

where the supremum is taken over all measurable sets. The probability distribution of a random variable $X$ is denoted by $\mathcal{L}(X)$ and that of Poisson distribution with parameter $\lambda$ by $\mathcal{P}(\lambda)$. Moreover, if $f(n)$ and $g(n)$ are two positive functions, notation $f(n) = \mathcal{O}(g(n))$ means that there exists a constant $C > 0$ such that, for all $n$, $f(n) \leq Cg(n)$.

Theorem 3.1 states that the asymptotic behavior of $X_n(\eta)$ is Poissonian and defines the speed of convergence in terms of total variation distance.

THEOREM 3.1. *Let $\eta$ be a local configuration of radius $r$. Assume that the magnetic field $a(n)$ is negative, the pair potential $b(n)$ is nonnegative and they satisfy*

$$n^d W_n(\eta) = \lambda,$$

14 D. COUPIER*for some constant $\lambda > 0$. Furthermore, if the hypotheses* (H1) *and* (H2) *are satisfied, then the total variation distance between the distribution of $X_n(\eta)$ and the Poisson distribution with parameter $\lambda$ is such that*:

$$d_{\mathrm{TV}}(\mathcal{L}(X_n(\eta)), \mathcal{P}(\lambda)) = \mathcal{O}(\max\{\Delta_n(\eta), \Theta_n(\eta)\}).$$

Now, some comments are needed. First, this result generalizes the Poisson approximation given in [6] (Theorem 1.3). When the pair potential $b(n) = b > 0$ is fixed, the product $n^d \exp(2a(n)k(\eta))$ becomes constant. Denote by $\lambda'$ this quantity. Then, the limit distribution for $X_n(\eta)$ is the Poisson distribution with parameter $\lambda = \lambda' \exp(-2b(n)\gamma(\eta))$. Furthermore, the probability gap $\Delta_n(\eta)$ and the maximality probability $\Theta_n(\eta)$ are of order $\exp(2a(n))$ and consequently,

$$d_{\mathrm{TV}}(\mathcal{L}(X_n(\eta)), \mathcal{P}(\lambda' \exp(-2b(n)\gamma(\eta)))) = \mathcal{O}(n^{-d/k(\eta)}).$$

We believe that $\max\{\Delta_n(\eta), \Theta_n(\eta)\}$ is the real speed at which the total variation distance between $\mathcal{L}(X_n(\eta))$ and $\mathcal{P}(\lambda)$ tends to zero. Indeed, it seems to be true for the upper bound given by Lemma 4.5 (for more details, see Chapter 3 of [2]). Moreover, this has been proved by Ganesh et al. [15] in the case where the local configuration $\eta$ represents a single positive vertex [with $k(\eta) = 1$, $\gamma(\eta) = \mathcal{V}$] and the pair potential is fixed.

Under the hypotheses (H1) and (H2), the expectation of $X_n(\eta)$ satisfies

$$n^d W_n(\eta)(1 - \mathcal{O}(\Delta_n(\eta))) \leq \mathbb{E}_{a,b}[X_n(\eta)] \leq n^d W_n(\eta)(1 + \mathcal{O}(\Delta_n(\eta)))$$

[see relations (20) and (21) of the next section]. In particular, if the product $n^d W_n(\eta)$ tends to 0 (resp. $+\infty$), the same is true for $\mathbb{E}_{a,b}[X_n(\eta)]$. Actually, a better result can be easily deduced from Theorem 3.1: if the product $n^d W_n(\eta)$ tends to 0 (resp. $+\infty$), then the probability $\mu_{a,b}(X_n(\eta) > 0)$ tends to 0 (resp. 1). Roughly speaking, if $W_n(\eta)$ is small compared to $n^{-d}$, then asymptotically, there is no occurrence of $\eta$ in $G_n$. If $W_n(\eta)$ is large compared to $n^{-d}$, then at least one occurrence of $\eta$ can be found in the graph, with probability tending to 1. Using the vocabulary of the random graph theory, this means that the quantity $n^{-d}$ is the *threshold function* for $W_n(\eta)$ of the property $X_n(\eta) > 0$. A straight proof of this statement can be found in [8]. A result similar to (actually weaker than) Theorem 3.1 can be proved without using the Stein–Chen method. Under the same hypotheses, the distribution of $X_n(\eta)$ converges weakly to $\mathcal{P}(\lambda)$, as $n$ tends to infinity. The proof of this result is based on the moment method (see [3], page 25) and Lemma 2.4; see the ideas developed in the proof of Lemma 4.6 for details.

Replacing Proposition 4.4 with the Stein–Chen formulation of [1], it seems the estimates given by Lemma 2.4 should apply to systems formed by "animals" interacting only by volume exclusions (as in [12] and [13]).

Finally, as has been suggested by an anonymous referee, Theorem 3.1 can be extended to some nonlocal configurations whose cardinalities diverge suitably with the size of the lattice.

The rest of this section is devoted to translating the hypotheses of Theorem 3.1 in terms of magnetic field $a(n)$ and pair potential $b(n)$.

PROPOSITION 3.2. *Let $\eta$ be a local configuration of radius $r$. Assume that the magnetic field $a(n)$ is negative, the pair potential $b(n)$ is nonnegative and they satisfy*

$$n^d W_n(\eta) = \lambda,$$

*for some constant $\lambda > 0$. Then the two conditions $a(n) + 2\mathcal{V}b(n) \leq 0$, for $n$ large enough, and*

$$\lim_{n \to +\infty} a(n) + \mathcal{V}b(n) = -\infty$$

*imply respectively the two hypotheses* (H2) *and* (H1).

The limit $a(n) + \mathcal{V}b(n) \to -\infty$ forces the magnetic field $a(n)$ to tend to $-\infty$. Proposition 3.2 and Theorem 3.1 state that $a(n) \to -\infty$ is the main condition bearing on the magnetic field $a(n)$ in order to expect a Poisson approximation. Rare positive vertices among a majority of negative ones constitutes the global context in which our Poisson approximations take place. Conditions bearing on the pair potential $b(n)$ are much larger. Here are two examples of couples $(a(n), b(n))$ satisfying the conditions of Proposition 3.2 and so the hypotheses of Theorem 3.1.

EXAMPLE 3.3. Assume the potentials $a(n) < 0$ and $b(n) \geq 0$ are such that $n^d W_n(\eta) = \lambda$, for some $\lambda > 0$, and

$$\lim_{n \to +\infty} a(n) = -\infty \quad \text{and} \quad \lim_{n \to +\infty} \frac{b(n)}{a(n)} = 0.$$

Then, conditions of Proposition 3.2 are satisfied. Indeed, the quantity

$$a(n) + \mathcal{V}b(n) = a(n)\left(1 + \mathcal{V}\frac{b(n)}{a(n)}\right)$$

tends to $-\infty$, so does $a(n) + 2\mathcal{V}b(n)$. The pair potential $b(n)$ is allowed to tend to $+\infty$ [but slower than $|a(n)|$] or to be bounded.

EXAMPLE 3.4. Assume the potentials $a(n)$ and $b(n)$ satisfy

$$a(n) = \frac{1}{2k(\eta) + \gamma(\eta)/\mathcal{V}} \log\left(\frac{\lambda}{n^d}\right),$$

$$b(n) = -\frac{a(n)}{2\mathcal{V}} = \frac{-1}{4\mathcal{V}k(\eta) + 2\gamma(\eta)} \log\left(\frac{\lambda}{n^d}\right),$$

where $\lambda$ is a positive constant. The rate at which the pair potential $b(n)$ tends to infinity is the same as $|a(n)|$. It is easy to check that the magnetic field $a(n)$ is negative, $b(n)$ is nonnegative, the product $n^d W_n(\eta)$ is equal to $\lambda$, the quantity $a(n) + 2\mathcal{V}b(n)$ is null and $a(n) + \mathcal{V}b(n) = a(n)/2$ tends to $-\infty$. In other words, conditions of Proposition 3.2 are satisfied.

Let us end this section by proving Proposition 3.2.

PROOF OF PROPOSITION 3.2.   Let $\eta$ be a local configuration of radius $r$. Assume that potentials $a(n)$ and $b(n)$ satisfy conditions of Proposition 3.2.

First, let us prove (H1) is satisfied. The study of the probability gap has been already done; $\Delta_n(\eta) \leq \exp(2a(n))$. Hence, the limit $a(n) + \mathcal{V}b(n) \to -\infty$ forces the magnetic field $a(n)$ to tend to $-\infty$ and $\Delta_n(\eta)$ to 0. Now, let us bound the maximality probability $\Theta_n(\eta)$. Let $\eta' \in \mathcal{C}_r^*(\eta)$. Although the set of positive vertices of $\eta'$ contains that of $\eta$, the perimeter $\gamma(\eta')$ is not necessarily larger than $\gamma(\eta)$: roughly, $V_+(\eta)$ may have holes. However, the inequality

$$\gamma(\eta') \geq \gamma(\eta) - (k(\eta') - k(\eta))\mathcal{V}$$

holds. Hence, the ratio $W_n(\eta')$ divided by $W_n(\eta)$ is bounded:

$$\frac{W_n(\eta')}{W_n(\eta)} = \exp(2a(n)(k(\eta') - k(\eta)) - 2b(n)(\gamma(\eta') - \gamma(\eta)))$$

$$\leq \exp((2a(n) + 2\mathcal{V}b(n))(k(\eta') - k(\eta))).$$

Then the limit $a(n) + \mathcal{V}b(n) \to -\infty$ implies the maximality gap satisfies

$$\Theta_n(\eta) = \max_{\eta' \in \mathcal{C}_r^*(\eta)} \frac{W_n(\eta')}{W_n(\eta)}$$

$$\leq \exp(2a(n) + 2\mathcal{V}b(n))$$

and tends to 0 as $n$ tends to infinity.

It remains to prove that (H2) holds. Let $V$ be a set of vertices and let $\zeta$ be a configuration of $\{-,+\}^V$. Let us start with writing

$$\mu_{a,b}(I_V^\zeta = 1) = \sum_{\sigma \in \{-,+\}^{\delta V}} \mu_{a,b}(I_V^\zeta = I_{\delta V}^\sigma = 1).$$

Let $\sigma$ and $\omega$ be two configurations respectively defined on $\delta V$ and $\delta \overline{V}$. Thanks to Lemma 2.2 and relation (12), we get

$$\mu_{a,b}(I_V^{\zeta\sigma} = 1 | I_{\delta \overline{V}}^\omega = 1) = \frac{W_n(\zeta\sigma\omega)}{\sum_{\vartheta \in \{-,+\}^{\overline{V}}} W_n(\vartheta\omega)}$$

$$\leq \frac{W_n(\zeta\sigma\omega)}{W_n(\omega)}$$

$$\leq W_n(\zeta\sigma)\exp(4b(n)\operatorname{conn}(\zeta\sigma,\omega))$$
$$\leq W_n(\zeta\sigma)\exp(4b(n)\operatorname{conn}(\sigma,\omega)).$$

Moreover, the weight of the configuration $\zeta\sigma$ can be expressed as

$$W_n(\zeta\sigma) = W_n(\zeta)\exp(2a(n)k(\sigma) - 2b(n)(\gamma(\zeta\sigma) - \gamma(\zeta)))$$
$$= W_n(\zeta)\exp(2a(n)k(\sigma) - 2b(n)(\gamma(\sigma) - 2\operatorname{conn}(\zeta,\sigma))).$$

Combining the previous relations, we deduce the following inequality:

$$\mu_{a,b}(I_V^{\zeta\sigma} = 1 | I_{\delta\overline{V}}^{\omega} = 1)$$
$$\leq W_n(\zeta)\exp(2a(n)k(\sigma) - 2b(n)(\gamma(\sigma) - 2\operatorname{conn}(\zeta,\sigma) - 2\operatorname{conn}(\sigma,\omega))).$$

A way to obtain an upper bound for the above conditional probability which does not depend on $\omega$ consists in using the inequality

$$\operatorname{conn}(\zeta,\sigma) + \operatorname{conn}(\sigma,\omega) \leq \mathcal{V}k(\sigma).$$

As a consequence,

$$\mu_{a,b}(I_V^{\zeta\sigma} = 1 | I_{\delta\overline{V}}^{\omega} = 1) \leq W_n(\zeta)\exp((2a(n) + 4\mathcal{V}b(n))k(\sigma) - 2b(n)\gamma(\sigma))$$
$$\leq W_n(\zeta),$$

whenever $a(n) + 2\mathcal{V}b(n)$ is negative. Then, the probability $\mu_{a,b}(I_V^{\zeta\sigma} = 1)$ is bounded by the weight $W_n(\zeta)$ and finally we get

$$\mu_{a,b}(I_V^{\zeta} = 1) \leq 2^{|\delta V|} W_n(\zeta).$$

Observe that the constant $C(V) = 2^{\mathcal{V}|V|} \geq 2^{|\delta V|}$ is suitable and depends on the set $V$ only through its cardinality. □

**4. Proof of Theorem 3.1.** This section is devoted to the proof of Theorem 3.1, whose layout is essentially the same as the one of the proof of Theorem 1.3 of [6] [i.e., in the case $b(n) = b > 0$]. Let $\eta$ be a local configuration of radius $r$. Throughout this section, we shall suppose the magnetic field $a(n)$ is negative, the pair potential $b(n)$ is nonnegative and they satisfy

$$n^d W_n(\eta) = \lambda,$$

for some constant $\lambda > 0$.

For all vertex $x$, the indicator random variable $\overline{I}_x^{\eta}$ is defined by

$$\overline{I}_x^{\eta} = \sum_{\eta' \in \mathcal{C}_r(\eta)} I_x^{\eta'},$$

where $\mathcal{C}_r(\eta)$ is formed by the local configurations of radius $r$ whose set of positive vertices contains that of $\eta$. Thus, we denote by $\overline{X}_n(\eta)$ the sum of these indicators over $x \in V_n$:

$$\overline{X}_n(\eta) = \sum_{x \in V_n} \overline{I}_x^{\eta} \tag{16}$$
$$= X_n(\eta) + \sum_{\eta' \in \mathcal{C}_r^*(\eta)} X_n(\eta'),$$

whose expectation $\mathbb{E}_{a,b}[\overline{X}_n(\eta)]$ will be simply denoted by $\lambda_n$. Finally, in order to simplify formulas, the maximum of the probability gap $\Delta_n(\eta)$ and the maximality probability $\Theta_n(\eta)$ will be merely denoted by $M_n(\eta)$:

$$M_n(\eta) = \max\{\Delta_n(\eta), \Theta_n(\eta)\}.$$

The proof of Theorem 3.1 is organized as follows. The total variation distance between $\mathcal{L}(X_n(\eta))$ and $\mathcal{P}(\lambda)$ is bounded by

$$d_{\mathrm{TV}}(\mathcal{L}(X_n(\eta)), \mathcal{L}(\overline{X}_n(\eta))) + d_{\mathrm{TV}}(\mathcal{L}(\overline{X}_n(\eta)), \mathcal{P}(\lambda_n)) + d_{\mathrm{TV}}(\mathcal{P}(\lambda_n), \mathcal{P}(\lambda)).$$

We are going to prove that each term of the above sum is of order $\mathcal{O}(M_n(\eta))$. The first and the last ones are respectively dealt with using Lemmas 4.1 and 4.3. Applied to the family of indicators $\{\overline{I}_x^{\eta}, x \in V_n\}$, the Stein–Chen method gives an upper bound for the second term (Lemma 4.5). Finally, Lemma 4.6 implies that this upper bound is a $\mathcal{O}(M_n(\eta))$.

As the maximality probability $\Theta_n(\eta)$ tends to 0, the occurrences of local configurations of $\mathcal{C}_r^*(\eta) = \mathcal{C}_r(\eta) \setminus \{\eta\}$ have vanishing probability. Hence, the random variables $X_n(\eta)$ and $\overline{X}_n(\eta)$ (resp. their expectations $\mathbb{E}[X_n(\eta)]$ and $\lambda_n$) will be asymptotically equal.

LEMMA 4.1. *The total variation distance between the distributions of $X_n(\eta)$ and $\overline{X}_n(\eta)$ satisfies*

$$d_{\mathrm{TV}}(\mathcal{L}(X_n(\eta)), \mathcal{L}(\overline{X}_n(\eta))) = \mathcal{O}(M_n(\eta)). \tag{17}$$

*Furthermore, the following inequalities hold:*

$$\mathbb{E}_{a,b}[X_n(\eta)] \leq \lambda_n \leq \mathbb{E}_{a,b}[X_n(\eta)] + |\mathcal{C}_r^*(\eta)|\lambda C(B)\Theta_n(\eta), \tag{18}$$

*where $C(B)$ denotes the constant of the hypothesis* (H2) *corresponding to a ball of radius $r$.*

PROOF. The total variation distance between two probability distributions $\mu$ and $\nu$ can be written as

$$d_{\mathrm{TV}}(\mu, \nu) = \inf\{\mathbb{P}(X \neq Y), \mathcal{L}(X) = \mu \text{ and } \mathcal{L}(Y) = \nu\}.$$

Using this characterization and the identity (16), it follows that

$$d_{\text{TV}}(\mathcal{L}(X_n(\eta)), \mathcal{L}(\overline{X}_n(\eta))) \leq \mu_{a,b}(X_n(\eta) \neq \overline{X}_n(\eta))$$
$$\leq \mu_{a,b}(\overline{X}_n(\eta) > X_n(\eta))$$
$$\leq \mu_{a,b}\left(\sum_{\eta' \in \mathcal{C}_r^*(\eta)} X_n(\eta') > 0\right).$$

The above sum is an integer-valued variable. So, its probability of being positive is bounded by its expectation. Hence,

$$d_{\text{TV}}(\mathcal{L}(X_n(\eta)), \mathcal{L}(\overline{X}_n(\eta))) \leq \sum_{\eta' \in \mathcal{C}_r^*(\eta)} \mathbb{E}_{a,b}[X_n(\eta')].$$

Let $\eta' \in \mathcal{C}_r^*(\eta)$. The hypothesis (H2) allows to control the expectation of $\eta'$ with the maximality probability $\Theta_n(\eta)$:

$$\mathbb{E}_{a,b}[X_n(\eta')] = n^d \mu_{a,b}(I_x^{\eta'} = 1)$$
$$\leq n^d C(B) W_n(\eta')$$
$$\leq \lambda C(B) \frac{W_n(\eta')}{W_n(\eta)}$$
$$\leq \lambda C(B) \Theta_n(\eta).$$

Finally, we get

$$d_{\text{TV}}(\mathcal{L}(X_n(\eta)), \mathcal{L}(\overline{X}_n(\eta))) \leq |\mathcal{C}_r^*(\eta)| \lambda C(B) \Theta_n(\eta)$$

and (17) follows. Relation (18) is an immediate consequence of the previous inequalities and the identity (16). □

Lemmas 2.4 and 4.1 state that the expectation of $\overline{X}_n(\eta)$, that is, $\lambda_n$, is close to that of $X_n(\eta)$, which is itself close to $\lambda$.

LEMMA 4.2. *There exists a constant $K(r) > 0$ such that*

(19) $$\lambda(1 - K(r) M_n(\eta)) \leq \lambda_n \leq \lambda(1 + K(r) M_n(\eta)).$$

The constant $K(r)$ depends on the radius $r$ and on the parameters of the model (i.e., $\rho$, $p$ and the dimension $d$).

PROOF OF LEMMA 4.2. Let $x \in V_n$ and denote by $B$ the ball $B(x, r)$. First, note that the indicator $\mathbb{1}_{\sigma \neq \sigma^0}$ occurring in Lemma 2.4 can be expressed as the following sum:

$$\sum_{\sigma \neq \sigma^0} I_{\delta B}^\sigma,$$

where $\sigma^0$ represents the null configuration of $\{-,+\}^{\delta B}$. So, the upper bound (15) implies

$$\mathbb{E}_{a,b}[X_n(\eta)] = \mathbb{E}_{a,b}[n^d \mu_{a,b}(I_x^\eta = 1 | \mathcal{F}(\delta B))]$$
$$\leq \lambda \left(1 + \Delta_n(\eta) \sum_{\sigma \neq \sigma^0} \frac{\mathbb{E}_{a,b}[I_{\delta B}^\sigma]}{W_n(\sigma)}\right).$$

Let $C(\delta B)$ be the constant of the hypothesis (H2) corresponding to the neighborhood of a ball of radius $r$. Then,

$$\mathbb{E}_{a,b}[I_{\delta B}^\sigma] \leq C(\delta B) W_n(\sigma).$$

Hence, we obtain an upper bound for $\mathbb{E}_{a,b}[X_n(\eta)]$:

(20) $$\mathbb{E}_{a,b}[X_n(\eta)] \leq \lambda(1 + \Delta_n(\eta) 2^{\delta B} C(\delta B)).$$

Let us deal with the lower bound. It has been already said (Lemma 2.3) that the weight of a configuration $\sigma$ different from $\sigma^0$ is smaller than the probability gap $\Delta_n(\eta)$. As a consequence,

$$\mathbb{E}_{a,b}\left[\sum_{\sigma \neq \sigma^0} I_{\delta B}^\sigma\right] \leq \sum_{\sigma \neq \sigma^0} C(\delta B) W_n(\sigma)$$
$$\leq 2^{\delta B} C(\delta B) \Delta_n(\eta).$$

Next, we deduce from the lower bound (14):

(21)
$$\mathbb{E}_{a,b}[X_n(\eta)] = \mathbb{E}_{a,b}[n^d \mu_{a,b}(I_x^\eta = 1 | \mathcal{F}(\delta B))]$$
$$\geq \lambda(1 - |\mathcal{C}_r| \Delta_n(\eta))\left(1 - \mathbb{E}_{a,b}\left[\sum_{\sigma \neq \sigma^0} I_{\delta B}^\sigma\right]\right)$$
$$\geq \lambda(1 - |\mathcal{C}_r| \Delta_n(\eta))(1 - 2^{\delta B} C(\delta B) \Delta_n(\eta))$$
$$\geq \lambda(1 - \Delta_n(\eta)(|\mathcal{C}_r| + 2^{\delta B} C(\delta B))).$$

Finally, combining the inequalities between the expectation of $X_n(\eta)$ and $\lambda_n$ given by relation (18) and the inequalities between the expectation of $X_n(\eta)$ and $\lambda$ given by (20) and (21), we conclude with

$$K(r) = \max\{|\mathcal{C}_r| + 2^{\delta B} C(\delta B), |\mathcal{C}_r^*(\eta)| C(B) + 2^{\delta B} C(\delta B)\}. \qquad \square$$

The total variation distance between two probability distributions on the set of integers can be expressed as

$$d_{\mathrm{TV}}(\mu, \nu) = \tfrac{1}{2} \sum_{m \geq 1} |\mu(m) - \nu(m)|.$$

Using this characterization, relation (19) between $\lambda_n$ and $\lambda$ and hypothesis (H1), it is possible to bound the total variation distance between $\mathcal{P}(\lambda_n)$ and $\mathcal{P}(\lambda)$. This result is very similar to Lemma 4.2 of [6]. So, it will not be proved here.

LEMMA 4.3. *The total variation distance between the Poisson distributions with parameters $\lambda_n$ and $\lambda$ satisfies*

$$d_{\mathrm{TV}}(\mathcal{P}(\lambda_n), \mathcal{P}(\lambda)) = \mathcal{O}(M_n(\eta)).$$

There remains to bound the total variation distance between $\mathcal{L}(\overline{X}_n(\eta))$ and $\mathcal{P}(\lambda_n)$. This is based on the Stein–Chen method and particularly on Corollary 2.C.4, page 26 of [2] which is described below (Proposition 4.4). Let $\{I_i\}_{i \in I}$ be a family of random indicators with expectations $\pi_i$. Let us denote

$$Z = \sum_{i \in I} I_i \quad \text{and} \quad \theta = \sum_{i \in I} \pi_i.$$

The random variables $\{I_i\}_{i \in I}$ are *positively related* if for each $i$, there exist random variables $\{J_{j,i}\}_{j \in I}$ defined on the same probability space such that

$$\mathcal{L}(J_{j,i}, j \in I) = \mathcal{L}(I_j, j \in I | I_i = 1)$$

and, for all $j \neq i$, $J_{j,i} \geq I_j$.

PROPOSITION 4.4. *If the random variables $\{I_i\}_{i \in I}$ are positively related, then*

$$d_{\mathrm{TV}}(\mathcal{L}(Z), \mathcal{P}(\theta)) \leq \frac{1 - e^{-\theta}}{\theta} \left( \mathrm{Var}(Z) - \theta + 2 \sum_{i \in I} \pi_i^2 \right).$$

Proposition 4.4 can be applied to our context. First, observe there is a natural partial ordering on the configuration set $\mathcal{X}_n = \{-, +\}^{V_n}$ defined by $\sigma \leq \sigma'$ if $\sigma(x) \leq \sigma'(x)$ for all vertices $x \in V_n$. A function $f : \mathcal{X}_n \to \mathbb{R}$ is *increasing* if $f(\sigma) \leq f(\sigma')$ whenever $\sigma \leq \sigma'$. By construction, the indicators $\overline{I}_x^\eta$ are increasing functions. Furthermore, for a positive value of the pair potential $b(n)$, the Gibbs measure $\mu_{a,b}$ defined by (2) satisfies the FKG inequality, that is,

(22) $$\mathbb{E}_{a,b}[fg] \geq \mathbb{E}_{a,b}[f] \mathbb{E}_{a,b}[g],$$

for all increasing functions $f$ and $g$ on $\mathcal{X}_n$; see, for instance, Section 3 of [14]. Then, Theorem 2.G, page 29 of [2] implies that the increasing random indicators $\overline{I}_x^\eta$, $x \in V_n$, are positively related. Replacing $I_i$ with $\overline{I}_x^\eta$, $Z$ with $\overline{X}_n(\eta)$ and $\theta$ with $\lambda_n$, Proposition 4.4 produces the following result.

LEMMA 4.5. *The following inequality holds:*

$$d_{\mathrm{TV}}(\mathcal{L}(\overline{X}_n(\eta)), \mathcal{P}(\lambda_n)) \tag{23}$$
$$\leq \frac{1}{\lambda_n}\left(\mathrm{Var}_{a,b}[\overline{X}_n(\eta)] - \lambda_n + 2\sum_{x \in V_n} \mathbb{E}_{a,b}[\overline{I}_x^\eta]^2\right).$$

The last term of the above bound is equal to the ratio $\lambda_n^2$ divided by $n^d$:

$$\sum_{x \in V_n} \mathbb{E}_{a,b}[\overline{I}_x^\eta]^2 = \frac{\lambda_n^2}{n^d}.$$

So, it tends to 0 as $n$ tends to infinity. Hence, the bound (23) indicates that, as $n \to +\infty$, the distance to the Poisson approximation is essentially the difference between the variance and the expectation of $\overline{X}_n(\eta)$. This difference can be written as $M_2(\overline{X}_n(\eta)) + \lambda_n^2$ where $M_2(\overline{X}_n(\eta))$ denotes the second moment of the random variable $\overline{X}_n(\eta)$:

$$M_2(\overline{X}_n(\eta)) = \mathbb{E}_{a,b}[\overline{X}_n(\eta)(\overline{X}_n(\eta) - 1)].$$

The quantity $\lambda_n^2$ is controlled by Lemma 4.2, and the second moment $\overline{X}_n(\eta)$ by the following result.

LEMMA 4.6. *The second moment of the random variable $\overline{X}_n(\eta)$ satisfies*

$$M_2(\overline{X}_n(\eta)) = \lambda^2 + \mathcal{O}(M_n(\eta)).$$

We deduce from Lemma 4.5 that

$$d_{\mathrm{TV}}(\mathcal{L}(\overline{X}_n(\eta)), \mathcal{P}(\lambda_n)) \leq \frac{1}{\lambda_n}\left(M_2(\overline{X}_n(\eta)) + \lambda_n^2\left(\frac{2}{n^d} - 1\right)\right).$$

Thus, combining inequalities given by Lemmas 4.2 and 4.6, it follows that

$$d_{\mathrm{TV}}(\mathcal{L}(\overline{X}_n(\eta)), \mathcal{P}(\lambda_n)) = \mathcal{O}(M_n(\eta) + n^{-d}).$$

We conclude using $\lambda n^{-d} = W_n(\eta) \leq \Delta_n(\eta) \leq M_n(\eta)$.

Let us finish the proof of Theorem 3.1 by proving Lemma 4.6.

PROOF OF LEMMA 4.5. The variable $\overline{X}_n(\eta)$ counts the number of copies in the graph $G_n$ of local configurations belonging to $\mathcal{C}_r(\eta)$. So, the quantity $M_2(\overline{X}_n(\eta))$ can be interpreted as the expected number of ordered couples of copies of elements of $\mathcal{C}_r(\eta)$:

$$M_2(\overline{X}_n(\eta)) = \sum_{k \geq 2} \mu_{a,b}(\overline{X}_n(\eta) = k)\frac{k!}{(k-2)!}$$

$$= \mathbb{E}_{a,b}\left[\sum_{(x_1, x_2) \in V_n^2} \overline{I}_{x_1}^\eta \times \overline{I}_{x_2}^\eta\right].$$

Let $(x_1, x_2)$ be a couple of vertices. From now on, two cases must be distinguished: either the distance between $x_1$ and $x_2$ is smaller than $2r+3$ or not. The meaning of this distinction will be revealed later. However, we can already remark that the number of couples $(x_1, x_2)$ such that $\text{dist}(x_1, x_2) \leq 2r+3$ is bounded by $\beta(2r+3)n^d$. Indeed, there are $n^d$ possibilities for the first vertex $x_1$ and no more than $\beta(2r+3)$ possibilities for the second one since it belongs to the ball $B(x_1, 2r+3)$. As for the number of couples $(x_1, x_2)$ such that $\text{dist}(x_1, x_2) > 2r+3$, it is merely bounded by $n^{2d}$.

Each indicator $\overline{I}_x^\eta$ is defined as the sum of $I_x^{\eta'}$, $\eta' \in \mathcal{C}_r(\eta)$. So, the second moment $M_2(\overline{X}_n(\eta))$ can be expressed as

$$\sum_{\substack{\eta_1, \eta_2 \\ \in \mathcal{C}_r(\eta)}} \left( \mathbb{E}_{a,b}\left[ \sum_{\substack{(x_1,x_2) \in V_n^2 \\ \text{dist}(x_1,x_2) \leq 2r+3}} I_{x_1}^{\eta_1} \times I_{x_2}^{\eta_2} \right] + \mathbb{E}_{a,b}\left[ \sum_{\substack{(x_1,x_2) \in V_n^2 \\ \text{dist}(x_1,x_2) > 2r+3}} I_{x_1}^{\eta_1} \times I_{x_2}^{\eta_2} \right] \right)$$

where the above expectations will be respectively denoted by $E^{\leq}(\eta_1, \eta_2)$ and $E^{>}(\eta_1, \eta_2)$. We are going to prove the three following statements from which Lemma 4.6 follows:

(24) $\quad \forall (\eta_1, \eta_2) \in \mathcal{C}_r(\eta)^2 \quad E^{\leq}(\eta_1, \eta_2) = \mathcal{O}(M_n(\eta));$

(25) $\quad \hspace{10em} E^{>}(\eta, \eta) = \lambda^2 + \mathcal{O}(M_n(\eta));$

(26) $\quad \forall (\eta_1, \eta_2) \in \mathcal{C}_r(\eta)^2 \setminus \{(\eta, \eta)\} \quad E^{>}(\eta_1, \eta_2) = \mathcal{O}(M_n(\eta)).$

Let $(x_1, x_2)$ be a couple of vertices such that $\text{dist}(x_1, x_2) \leq 2r+3$ and let $\eta_1, \eta_2$ be two local configurations of $\mathcal{C}_r(\eta)$. The balls $B(x_1, r)$ and $B(x_2, r)$ are both included in the large ball $B(x_1, 3r+3)$. The intuition is that if $\eta_1$ and $\eta_2$ occur on $B(x_1, r)$ and $B(x_2, r)$, then locally [i.e., in the ball $B(x_1, 3r+3)$] at least $k(\eta) + 1$ positive vertices are present. This has a vanishing probability. So as to lighten formulas, let us denote by $B$ the ball $B(x_1, r)$ and by $R$ the set of vertices $y$ such that $r < \text{dist}(x_1, y) \leq 3r+3$. The event $I_{x_1}^{\eta_1} = I_{x_2}^{\eta_2} = 1$ implies

$$\left( \bigcup_{\eta' \in \mathcal{C}_r^*(\eta)} \bigcup_{\sigma \in \{-,+\}^R} I_B^{\eta'} = I_R^{\sigma} = 1 \right) \cup \left( \bigcup_{\sigma \in \{-,+\}^R \setminus \{\sigma^0\}} I_B^{\eta} = I_R^{\sigma} = 1 \right),$$

where $\sigma^0$ represents the null configuration of $\{-, +\}^R$. First, remark that upper bound of Lemma 2.4 can be extended from $\{-, +\}^{\delta B}$ to $\{-, +\}^R$. For all $\eta' \in \mathcal{C}_r(\eta)$ and $\sigma \in \{-, +\}^R$,

(27) $\quad \mu_{a,b}(I_B^{\eta'} = 1 | \sigma) \leq W_n(\eta') \left( \mathbb{1}_{\sigma = \sigma^0} + \mathbb{1}_{\sigma \neq \sigma^0} \frac{\Delta_n(\eta')}{W_n(\sigma)} \right).$

Let us consider $\eta' \in \mathcal{C}_r^*(\eta)$ and $\sigma \in \{-,+\}^R$. Relation (27) implies

$$\mu_{a,b}(I_B^{\eta'} = I_R^\sigma = 1) = \mu_{a,b}(I_B^{\eta'} = 1 | I_R^\sigma = 1) \mu_{a,b}(I_R^\sigma = 1)$$
$$\leq W_n(\eta')\left(1 + \frac{\Delta_n(\eta')}{W_n(\sigma)}\right)\mu_{a,b}(I_R^\sigma = 1)$$
$$\leq W_n(\eta')(1 + C(R)\Delta_n(\eta')),$$

where $C(R)$ is the constant of (H2) associated to the set $R$. Thus, using

$$\Delta_n(\eta') \leq \exp(2a(n)) \leq 1$$

and $\eta' \in \mathcal{C}_r^*(\eta)$, it follows that

$$\mu_{a,b}(I_B^{\eta'} = I_R^\sigma = 1) \leq W_n(\eta)(1 + C(R))\Theta_n(\eta).$$

Now, if the configuration $\sigma$ belonging to $\{-,+\}^R$ is different from $\sigma^0$, then

$$\mu_{a,b}(I_B^\eta = I_R^\sigma = 1) = \mu_{a,b}(I_B^\eta = 1 | I_R^\sigma = 1)\mu_{a,b}(I_R^\sigma = 1)$$
$$\leq W_n(\eta)\frac{\Delta_n(\eta)}{W_n(\sigma)}\mu_{a,b}(I_R^\sigma = 1)$$
$$\leq W_n(\eta)C(R)\Delta_n(\eta).$$

In conclusion, we obtain an explicit bound for $E^\leq(\eta_1, \eta_2)$:

$$E^\leq(\eta_1, \eta_2) \leq \beta(2r+3)\lambda(|\mathcal{C}_r^*(\eta)|2^{|R|}(1 + C(R))\Theta_n(\eta) + 2^{|R|}C(R)\Delta_n(\eta)),$$

from which relation (24) follows.

Let $(x_1, x_2)$ be a couple of vertices such that $\text{dist}(x_1, x_2) > 2r + 3$ and let $\eta_1, \eta_2$ be two local configurations of $\mathcal{C}_r(\eta)$. Here, denote respectively by $B_1$ and $B_2$ the balls $B(x_1, r)$ and $B(x_2, r)$. Since the distance between $x_1$ and $x_2$ is larger than $2r + 3$, no vertex of $B_1$ can be a neighbor of a vertex of $B_2$ [actually, $\text{dist}(x_1, x_2) > 2r + 1$ suffices]. The Gibbs measure $\mu_{a,b}$ yields a Markov random field with respect to neighborhoods defined in (1) (see, e.g., [17], Lemma 3, page 7). As a consequence,

$$\mu_{a,b}(I_{x_1}^{\eta_1} = I_{x_2}^{\eta_2} = 1) = \mathbb{E}_{a,b}[\mu_{a,b}(I_{x_1}^{\eta_1} = I_{x_2}^{\eta_2} = 1 | \mathcal{F}(\delta B_1 \cup \delta B_2))]$$
(28)
$$= \mathbb{E}_{a,b}\left[\prod_{i=1}^2 \mu_{a,b}(I_{x_i}^{\eta_i} = 1 | \mathcal{F}(\delta B_1 \cup \delta B_2))\right]$$
$$= \mathbb{E}_{a,b}\left[\prod_{i=1}^2 \mu_{a,b}(I_{x_i}^{\eta_i} = 1 | \mathcal{F}(\delta B_i))\right].$$

In a first time, assume that $\eta_1 = \eta_2 = \eta$. Let $\sigma$ be a configuration of $\mathcal{X}_n$. Thanks to Lemma 2.4, we can write

$$\prod_{i=1}^{2} \mu_{a,b}(I_{x_i}^{\eta} = 1 | \sigma_{\delta B_i})$$
$$\leq W_n(\eta)^2 \left(1 + \Delta_n(\eta)\left(\frac{1}{W_n(\sigma_{\delta B_1})} + \frac{1}{W_n(\sigma_{\delta B_2})}\right)\right.$$
$$\left. + \Delta_n(\eta)^2 \cdot \left(\frac{1}{W_n(\sigma_{\delta B_1}) W_n(\sigma_{\delta B_2})}\right)\right).$$

The reason for which the distance between vertices $x_1$ and $x_2$ is assumed larger than $2r + 3$ is the following: no vertex of $\delta B_1$ can be a neighbor of a vertex of $\delta B_2$. So, the product $W_n(\sigma_{\delta B_1}) W_n(\sigma_{\delta B_2})$ is equal to the weight $W_n(\sigma_{\delta B_1 \cup \delta B_2})$. Hence, the inequality $\Delta_n(\eta) \leq 1$ and relation (28) imply that the probability $\mu_{a,b}(I_{x_1}^{\eta} = I_{x_2}^{\eta} = 1)$ is bounded by

$$W_n(\eta)^2 \left(1 + \Delta_n(\eta) \sum_{\sigma \in \mathcal{X}_n} \left(\frac{1}{W_n(\sigma_{\delta B_1})} + \frac{1}{W_n(\sigma_{\delta B_2})} + \frac{1}{W_n(\sigma_{\delta B_1 \cup \delta B_2})}\right) \mu_{a,b}(\sigma)\right).$$

Let $C(\delta B)$ and $C(2\delta B)$ be the constants of (H2) respectively associated to the sets $\delta B_1$ (or $\delta B_2$) and $\delta B_1 \cup \delta B_2$. Then, a new upper bound for the quantity $\mu_{a,b}(I_{x_1}^{\eta} = I_{x_2}^{\eta} = 1)$ is obtained:

$$W_n(\eta)^2 (1 + \Delta_n(\eta)(2^{|\delta B_1|} C(\delta B) + 2^{|\delta B_2|} C(\delta B) + 2^{|\delta B_1| + |\delta B_2|} C(2\delta B))).$$

Finally, the term $E^>(\eta, \eta)$ satisfies the inequality

$$E^>(\eta, \eta) \leq n^{2d} W_n(\eta)^2 (1 + \mathcal{O}(\Delta_n(\eta)))$$

and (25) is proved.

Assume at least one of the two local configurations $\eta_1, \eta_2 \in \mathcal{C}_r(\eta)$ is different from $\eta$. It remains to prove that $E^>(\eta_1, \eta_2)$ is of order $M_n(\eta)$. Techniques used in the previous case allow us to write

$$\prod_{i=1}^{2} \mu_{a,b}(I_{x_i}^{\eta_i} = 1 | \mathcal{F}(\delta B_i)) = \mathcal{O}(W_n(\eta_1) W_n(\eta_2)).$$

Since this upper bound does not depend on the configuration on $\delta B_1 \cup \delta B_2$, we deduce from (28):

$$\mu_{a,b}(I_{x_1}^{\eta_1} = I_{x_2}^{\eta_2} = 1) = \mathcal{O}(W_n(\eta_1) W_n(\eta_2)).$$

Now, if the local configuration $\eta_i$, for $i = 1, 2$, is different from $\eta$, then its weight $W_n(\eta_i)$ is smaller than $W_n(\eta) \Theta_n(\eta)$. So, $\Theta_n(\eta) \leq 1$ for $n$ large enough and $(\eta_1, \eta_2) \neq (\eta, \eta)$ imply that the product $W_n(\eta_1) W_n(\eta_2)$ is bounded by $W_n(\eta)^2 \Theta_n(\eta)$ and consequently (26) is proved. □

**Acknowledgments.** I would like to thank Roberto Fernández and an anonymous referee for their valuable advice.

UFR de Mathématiques  
Laboratoire Paul Painlevé  
Université Lille 1  
Cité Scientifique—Bât. M2  
59655 Villeneuve d'Ascq Cedex  
France  
E-mail: david.coupier@math.univ-lille1.fr